\documentstyle{amsppt}
\pagewidth{125mm}
\pageheight{190mm}
\parindent=8mm
\frenchspacing \tenpoint
\TagsAsMath

\topmatter
\title 
Ladders on Fano Varieties
\endtitle
\author Florin Ambro
\endauthor

\rightheadtext{Good divisors}

\address Department of Mathematics, The Johns Hopkins University, 
3400 N. Charles Street, Baltimore MD 21218, US \endaddress
\email ambro\@chow.mat.jhu.edu \endemail


\dedicatory
Dedicated to Professor V. A. Iskovskikh on his 60th birthday
\enddedicatory

\abstract
We prove the existence of ladders on log Fano varieties
of coindex less than $4$, as an application of
adjunction and nonvanishing.
\endabstract

\endtopmatter

\document

\head 1. Introduction
\endhead

This paper is related to the following 
\proclaim{Good Divisor Problem} Construct a 
regular ladder for a nef and big Cartier divisor $H$ 
on a normal variety $X$. This means to find a ``good'' element 
$S \in |H|$ and then repeat for $(S,H|_S)$. Note
that $H|_S$ is still nef and big.
\endproclaim
A ``good'' member means here an irreducible reduced
normal variety having singularities close to its 
ambient space. For example, they are both in the same 
class of singularities (log terminal, canonical,
terminal, $\epsilon$-terminal, etc).
In this paper we find sufficient conditions for
the existence of such divisors in the class
of Kawamata log terminal singularities (klt). The main result is 
(see Section 2 for definitions and notations):

\proclaim{Main Theorem}
 Let $H$ be a nef and big Cartier divisor on a normal variety $X$ of
 dimension $n$. Assume there is a boundary $B_X$ on $X$ such that

$(X,B_X)$ is klt,

$-(K_X+B_X) \equiv (n-r+1)H, n-r+1>0$,

$r<4$.  

 Then $dim|H| \geq n-1$ and $|H|$ does not have fixed components.
 Moreover, $(X,B_X+S)$ is a pure log terminal pair
  for the general $S \in |H|$.
 
 In particular, $(S,B_S=B_X|_S)$ is a weak log Fano with $co(H|_S)=co(H)=r$ 
 if $n>r$, a log CY if $n=r$ and a weak log general type if 
 $n<r<n+1$.
\endproclaim

A weak log Fano (CY, or general type) variety is a log pair 
$(X,B_X)$ with klt 
singularites such that $-(K_X+B_X)$ is nef and big ($K_X+B_X \equiv 0$,
or $K_X+B_X$ is nef and big). 
Note that the second condition of the theorem is equivalent to
$$
H \equiv K_X+B_X +tH, t=n-r+2>1
$$ 
Iterating the process, $H$ gives rise to a  
regular pure log terminal (plt) ladder
$$(X,H)=(X_n,H_n)>(X_{n-1},H_{n-1})> \cdots >
(X_c,H_c),\ \ \ c=\lfloor r \rfloor -1
$$
i.e. a regular ladder in the usual sense 
such that there are boundaries $B_i$ on $X_i$ with $(X_i,B_i+X_{i-1})$ plt,
$B_{i-1}=B_i |_{X_{i-1}}$ and $-(K_{X_j}+B_j) \equiv (j-r+1)H_i$.

All members, except the last one, are weak Fanos and the coindex of 
$H_j$ is constant on the ladder.

Note that the coindex of the weak 
Fanos on the ladder might decrease, which is actually better, since
the geometry of Fanos is simpler for smaller coindex. For example the
only Fano variety of coindex $0$ is the projective space. 
Moreover, this Fano ladders are maximal: the last member is non Fano.

\proclaim{Corollary} Let $H$ be a Cartier divisor of coindex $r<4$ on a 
weak log
Fano $(X,B_X)$ and $S \in |H|$ a general element. 

1) If $0<r<3$, then $S$ is smooth in base points. 
  The base locus is at most a finite set of points.

2) If $3 \leq r <4$, then $dimBsl|H| \leq 1$, $S$ is smooth
  in the base curves (if any) and it has only isolated singularities in the
  base locus. 
 Moreover, if $X$ is smooth in $Bsl|H|$, then $S$ is smooth in $Bsl|H|$ for 
$n \geq 4$ and it might have (exactly) canonical double points for $n=3$.
\endproclaim

The main techniques are the
Adjunction Conjecture, Kawamata's log canonical centers and the 
relation between dimension of log centers and discrepancies (see 
~\cite{He},~\cite{Me1} and also Lemma 5.3).
The main theorem follows from the 
Adjunction Conjecture and the effective nonvanishing on weak Fanos and
log-surfaces we prove in Section 3 and 4. 
In a previous draft we extended a trick from ~\cite{Me} to reduce the problem 
to the cases were adjunction was already proven. Kawamata's new Theorem
2.4\cite{Ka3} shows actually that the good divisor problem
is implied by adjunction and nonvanishing.
\newline
 As a corollary, we obtain a unified treatment of some results
of  Alexeev ~\cite{Al},
Prokhorov ~\cite{Pro}, Mella ~\cite{Me1, Me2}, Reid ~\cite{Re} and 
Shokurov ~\cite{Sh1}. I would like to mention Mella's ~\cite{Me1}
and Prokhorov's ~\cite{Pro} papers, which I found very inspiring. 

I would like to thank Professor Kawamata for updating me with his latest
result on the Adjunction Conjecture. I am grateful to Professor 
Shokurov for his guidance and
his valuable suggestions, and for allowing me to include his
unpublished theorem in the Appendix. The notion of log canonical 
singularity (LCS) subscheme is also from his lecture notes.

\head 2. Background
\endhead

 The ground field is assumed of characteristic $0$. 
A {\sl variety} $X$ 
means an irreducible, reduced, normal scheme of finite type over $k$.
A {\sl (sub)boundary} $B_X= \sum d_i E_i$ on $X$ is a 
$\Bbb Q$-Weil divisor with 
($d_i \leq 1, \forall i$) $0 \leq d_i \leq 1, \forall i$. A
(sub)boundary is called strict if $d_i<1,
\forall i$.
A {\sl log pair} $(X,B_X)$ is a variety $X$ equipped
with a boundary $B_X$ such that $K_X+B_X$ is a 
$\Bbb Q$-Cartier divisor on $X$.

\subhead 2.1. Shokurov's log canonical singularity subscheme
\endsubhead

 Fix a log pair $(X,B_X)$.
For a log resolution $f:Y\to X$, let $B^Y$ on $Y$ such that 
$K_Y+B^Y\sim_{\Bbb Q} f^*(K_X+B_X)$. Define 
$${\Cal I}(X,B_X)=f_* {\Cal O}_Y (\lceil-B^Y \rceil).$$
Then ${\Cal I}(X,B_X)$ is an ideal sheaf on $X$, independent of the 
choice of log resolution. The associated subscheme of $X$, denoted 
${\Cal L}(X,B_X)$ is called the {\sl log canonical singularity (LCS) 
subscheme associated to} $(X,B_X)$.
We also define the 
{\sl locus of log canonical singularities of} $(X,B_X)$ to be
$LCS(X,B_X)=Supp {\Cal L}(X,B_X)$.

\subhead 2.2. Singularities
\endsubhead

{\bf lc} The pair $(X,B_X)$ is called {\sl log canonical} (denoted lc) 
     if $B^Y$ is a subboundary for any resolution $Y \to X$. This 
     implies that ${\Cal L}(X,B_X)$ is a reduced scheme. 

{\bf klt} The pair $(X,B_X)$ is called {\sl Kawamata log terminal} 
       (denoted klt) if $B^Y$ is a strict suboundary for any
     resolution.

{\bf plt} The pair $(X,B_X)$ is called {\sl pure log terminal}
      (denoted plt) if it is log canonical and Kawamata log terminal
      in codimension bigger than $1$ (this is equivalent to the fact
      that $B^Y$ is a subboundary with all reduced components
      nonexceptional, for any resolution $f: Y \to X$).
\proclaim{Definition}
If $(X,B_X)$ is log canonical and $D$ is a $\Bbb Q$-Cartier effective divisor, 
we define the {\sl log canonical threshold of $D$}, denoted 
$lct(D;X,B_X)$, to be the maximal $\gamma \geq 0$ such that 
$(X,B_X+\gamma D)$ is log canonical. Note that $lct(D;X,B_X) \le 1$ 
if $D$ is a nonzero integer divisor.
\endproclaim

\subhead 2.3. Log canonical centers \cite{Ka1}
\endsubhead

We call any irreducible component $E$ of 
$\lceil -B^Y \rceil$ a {\sl log canonical place (lc place)} and $center_X(E) 
=f(E)$ is called a {\sl center
of log canonicity (lc center)}. 
We denote by $CLC(X,B_X)$ the set of all centers of log canonicity.
\proclaim{Proposition}~\cite{Ka1}
Let $(X,B_X^o)$ be a klt pair and $B_X > B_X^o$ an effective 
$\Bbb Q$-Cartier divisor such that $(X,B_X)$ is log canonical. 
Then

i) $CLC(X,B_X)$ is a finite set closed under intersection:
 \newline if $W_1,W_2 \in CLC(X,B_X)$ then the same holds for 
any irreducible component of $W_1 \cap W_2.$

ii) There are minimal elements of $CLC(X,B_X)$, called 
     {\sl minimal log canonical centers}. They are normal varieties.
\endproclaim

\proclaim{Remark}
 Let $(X,B_X)$ be a klt pair, $H$ be a nef and big Cartier
 divisor on $X$, and $D \equiv \gamma H$ be a $\Bbb Q$-Cartier divisor 
 such that $(X,B_X+D)$ is log canonical with $Z \in CLC(X,B_X+D)$ a 
 minimal lc center. 
Then, for any $0<\epsilon \ll 1$, $D$ can be perturbed to a 
$\Bbb Q$-divisor $D'$ such that $D' \equiv (\gamma+\epsilon)H$, 
$CLC(X,B_X+D')=\{Z\}$ and there is just one
log canonical place over $Z$ (for a given resolution).

Indeed, this can be obtained taking $D'=(1-\epsilon)D+\delta a M$, where 
$M \sim_{\Bbb Q} H$.
\endproclaim

\newpage

\subhead 2.4. Adjunction 
\endsubhead
\proclaim{The Adjunction Conjecture}
Let $(X,B_X^o)$ be a klt pair and $B_X > B_X^o$ an effective 
$\Bbb Q$-Cartier divisor such that $(X,B_X)$ is log canonical.
Then given a minimal log canonical center $W$ in 
$CLC(X,B_X)$, there 
is an effective $\Bbb Q$-Cartier divisor $B_W$ on $W$ such that 

i) $(W,B_W)$ is klt,

ii)$(K_X +B_X)|_W \equiv K_W + B_W$.
\endproclaim

\proclaim{Remark}
 This is a weak form of the Adjunction Theorem 
which was proved in several cases by Kawamata 
(~\cite{Ka1},~\cite{Ka2}): $codW\leq 2$ or $dimW=2$. 
\endproclaim
The latest result, which after a small perturbation is
the above conjecture, can substitute the adjunction formula
in this paper (see Lemma 5.3).  
\proclaim{Theorem 1}~\cite{Ka3}
Let $(X,B_X^o)$ be a klt pair and $B_X > B_X^o$ an effective 
$\Bbb Q$-Cartier divisor such that $(X,B_X)$ is log canonical.
Let $W$ be a minimal center of log canonical singularities for
$(X,B_X)$. Let $H$ be an ample Cartier divisor on $X$ and 
$\epsilon$ a rational positive number. Then there 
is an effective $\Bbb Q$-Cartier divisor $B_W$ on $W$ such that 

i) $(W,B_W)$ is klt,

ii)$(K_X +B_X+\epsilon H)|_W \sim_{\Bbb Q} K_W + B_W$.
\newline In particular, $W$ has only rational singularities.
\endproclaim

\subhead 2.5. Ladders \cite{Fu}
\endsubhead
Let $V$ be an $n$-fold polarized by an 
ample line bundle $L$. A sequence 
$$(V,L)=(V_n,L_n)>(V_{n-1},L_{n-1})> \cdots >(V_c,L_c)$$
is called a {\sl ladder} if each $V_{j-1}$ ($j=c+1,\ldots,n$) is an 
irreducible and reduced 
member of $|L_j|$, where $L_j$ is the restriction of $L$ to $V_j$. 
The ladder is called {\sl regular} if each restriction map 
$r:H^0(V_j,L_j) \to 
 H^0(V_{j-1},L_{j-1})$ is surjective.
\newline The existence of the above regular ladder implies that 
the base locus of $|L|$ will coincide with the base locus of $|L_c|$, 
hence it has dimension at most $c-1$.

\subhead 2.6.(Weak) Fano varieties
\endsubhead

A (weak) log Fano variety is a klt pair $(X,B_X)$  
such that $-(K_X+B_X)$ is (nef and big) ample.
In particular, $Pic(X)$ is torsion free, so any Cartier divisor
$H$ with $-(K_X+B_X) \equiv tH, t>0$ is uniquely (up to linear
equivalence) determined by $t$. We call $t$ the {\sl numerical index}
 of $H$ and denote it $index=ind(H)$. 
It is known that $0<t\leq n+1$, where $n=dimX$ 
(see ~\cite{Sh2}, for example). 
It turns out that a natural invariant for 
ladders will be the {\sl numerical coindex of} $H$, defined as
$co(H)=n+1-ind(H)$. The {\sl index (coindex) of}  $X$ is defined as the 
minimum
(maximum) index (coindex) over all Cartier divisors 
numerically proportional with $-(K_X+B_X)$.
Note that $0 \le co(H) <n+1$.

\proclaim{Remark} It is possible to extend general 
    properties from log Fano to weak log Fano varieties.
 Indeed, given a weak log Fano $(X,B_X)$ and $H$ Cartier of 
coindex $r$, there is a birational contraction 
$\mu:(X,B_X) \to (Y,B_Y)$ and $H_Y$ ample Cartier
on $Y$ such that

a) $H=\mu^*(H_Y)$,

b) $(Y,B_Y)$ is a log Fano and $H_Y$ has coindex $r$,

c) $K_X+B_X=\mu^*(K_Y+B_Y)$.
\newline This is a consequence of the 
Contraction Theorem ~\cite{KMM}. Note also that
if $S \in |H|, S_Y \in |H_Y|$ with $S=\mu^*S_Y$, then 
$lct(S;X,B_X)=lct(S_Y;Y,B_Y)$ and the log canonical centers on 
$Y$ are images of log canonical centers on $X$. 
Moreover, $Bsl|H|=\mu^{-1}(Bsl|H_Y|)$.
\endproclaim

\head 3. Nonvanishing on Fano varieties
\endhead

 In this section, $(X,B_X)$ is a weak log Fano variety of dimension
$n$ and $H$ is a nef and big Cartier divisor such that 
$-(K_X+B_X) \equiv ind(H) H,ind(H)>0$. We will compute
$h^0(H)=dim H^0(X,{\Cal O}_X(H))$ in some cases.
\proclaim{Notations}
 
$r=co(H)=n-ind(H)+1,d=H^n \geq 1, \delta = B_X . H^{n-1} \geq 0$

$p(t)=\chi({\Cal O}_X(tH))=a_nt^n+a_{n-1}t^{n-1}+\cdots + a_0$

$a_n=\frac{d}{n!}$, $a_{n-1}=\frac{-K_X . H^{n-1}}{2\ (n-1)!}=
\frac{ind(H) d+ \delta}{2\ (n-1)!}$.

$p_j = (-1)^n \chi(-j H)=h^0(K_X+ j H)\ge 0$ (from Serre duality and
Kawamata-Viehweg vanishing)
\endproclaim

Since $tH \equiv (K_X+B_X)+(t+n-r+1)H$, Kawamata-Viehweg 
vanishing \cite{KMM} implies
$p(t)=h^0(tH)$ for $t \ge -(n-\lfloor r\rfloor)$. In particular, any
integer $-(n-\lfloor r\rfloor) \le t <0$ is a zero of $p$, and 
$p(0)=\chi({\Cal O}_X)=1$.
\newline Since $(X,B_X)$ is klt, $X$ has rational singularities 
\cite{KMM}, so Serre duality applies.
\proclaim{Lemma 1} If $n \leq 2$, then $dim|H| >0$.

i) The only log Fano curve is ${\Bbb P}^1$ (with boundary), and
                      $dim|H|=d$.

ii) On a log Del Pezzo we have $dim|H|=\frac{(4-r)d+ \delta}{2}$.
\endproclaim
\demo{Proof}
 The first part is obvious. In the case of a log Del Pezzo, we know the  
Riemann Roch formula
$$p(t)=(d/2)t^2+(-K_X H/2) t + 1.$$ 
and $h^0(H)=p(1)$, so we are done.
\hfill $\square$
\enddemo
\proclaim{Lemma 2} Let $n\geq 3$. Then $dim|H|>0$ if $r=co(H)<4$.
More precisely,
 
i) If $\lfloor r\rfloor \leq 2$, then $h^0(H)=n-1+
           \frac{d(4-r)+\delta}{2}$.

ii)If $\lfloor r\rfloor = 3$, then $h^0(H)=n-1
           + p_{n-2} + \frac{d(4-r)+\delta}{2}$.

iii) If $\lfloor r\rfloor = 4$, then
        $h^0(H)=(1-p_{n-3})(n-1)+ p_{n-2} + \frac{d(4-r)+\delta}{2}.$

\endproclaim
\proclaim{Remark}
Note that $d(4-r)+\delta=-(K_X+(n-3)H. H^{n-1})$. 
In particular, if $d(4-r)+\delta> 0$, then $p_{n-3}=0$. 
\endproclaim
\demo{Proof}

i)(from ~\cite{Al}). We have
$$p(t)= 
  \frac{d}{n !} (t+1) \cdots (t+n-2) (t^2+a t +\frac{n(n-1)}{d})$$
$$ =\frac{d}{n!}t^n + \frac{d}{n!}(a+ \frac{(n-2)(n-1)}{2})t^{n-1}+\ldots$$
Therefore $a=-1+\frac{n(d(4-r)+\delta)}{2d}$ and 
 $h^0(H)= n-1+ \frac{d(4-r) + \delta}{2}$.

ii)(cf. ~\cite{Pro} for $n=r=3$). We have
$$
\align 
&p(t)=\frac{d}{n !} (t+1) \cdots (t+n-3) (t^3+a t^2 + bt + 
    \frac{n(n-1)(n-2)}{d})\\
&=\frac{d}{n!}t^n + \frac{d}{n!}(a+ \frac{(n-3)(n-2)}{2})t^{n-1}+ \ldots
\endalign
$$
Therefore,
$$
a=\frac{(6-r)n}{2} -3 + \frac{\delta n}{2 d}.
$$
Now, 
$$
p_{n-2}
=(-1)^n p(-n+2)\
=-1+\frac{d}{n(n-1)} ((n-2)^2-a(n-2)+b),
$$
therefore
$$
\align
&h^0(H)=n-2+ \frac{d}{n(n-1)} (1+(n-1)a -(n-2)^2) + p_{n-2} +1
\\
&=n-1 + p_{n-2} + \frac{d(4-r)+\delta}{2}.
\endalign
$$
iii) We have 
$$
\align
&p(t)=\frac{d}{n !} (t+1) \cdots (t+n-4) (t^4+a t^3 + b t^2 + c t+
  \frac{n(n-1)(n-2)(n-3)}{d})
\\
&=\frac{d}{n!}t^n + \frac{d}{n!}(a+ \frac{(n-4)(n-3)}{2})t^{n-1}+\ldots\\
\endalign
$$
Therefore $a=-6+\frac{n(d(8-r)+\delta)}{2d}$. Note that
$$
\align
&p(1)=n-3+\frac{d}{n(n-1)(n-2)}(1+a+b+c)\\
&p_{n-3}=1+\frac{d}{n(n-1)(n-2)}[(n-3)^3-a(n-3)^2+b(n-3)-c]\\
&p_{n-2}=n-3+\frac{d}{n(n-1)}[(n-2)^3-a(n-2)^2+b(n-2)-c]
\endalign
$$
We use the identity $(n-2)[b(n-2)-c] -(n-1)[b(n-3)-c]=b+c,$ 
and get
$$
\align
&b+c= (n-1)(n-3)^3-(n-2)^4-a[(n-1)(n-3)^2-(n-2)^3]\\
&-\frac{n(n-1)(n-2)}{d}[ (n-1)(p_{n-3} -1) - p_{n-2} +n-3 ],
\endalign
$$
hence
$$
\align
&p(1)=(n-1)(1-p_{n-3})+p_{n-2}+\frac{d}{n}[-2(n-3)+a]\\
&= (n-1)(1-p_{n-3})+p_{n-2} + \frac{d(4-r)+\delta}{2}.
\endalign
$$
\hfill $\square$
\enddemo
\proclaim{Remark}
We don't know if $dim|H|>0$, or even if $h^0(H)>0$ for all 
coindex $4$ divisors
$H$ on a (weak) log Fano (hence $n \geq 4$).

In this case, $h^0(H)=(n-1)(1-p_{n-3})+p_{n-2}+\frac{\delta}{2}$, and since
$-\delta=(K_X+(n-3)H.H^{n-1})$, we have $p_{n-3}=0$ unless
$\delta=B_X.H^{n-1}=0$. Therefore 

a) On a weak log Fano $n$-fold, a nef and big Cartier divisor
    $H$ of coindex $4$ such that $B_X.H^{n-1}>0$ has $dim|H| \geq n-2$. 

b) On a log Fano $n$-fold, an ample $H$
of coindex $4$ has $dim|H| \geq n-2$, unless $-K_X \sim (n-3)H, B_X=0$.

Let's assume $X$ is an exception as in b) above. From duality, 
$p(t)=p(-t-n+3)$, 
which is equivalent to $n$ even and there is $\alpha \in \Bbb Q$ 
such that the residual polynomial of degree $4$ appearing in the 
expression of $p(t)$ has the form
$$
t^4+2(n-3) t^3+\alpha t^2+[(n-3)\alpha -(n-3)^3]t+\frac{n(n-1)(n-2)(n-3)}{d}.
$$
Assuming now that $X$ is a $4$-fold , we get
$$
p(t)=p_k(t)=\frac{d}{4!} t(t+1)(t-1)(t+2)+\frac{k}{2}t(t+1)+1, 
 k \in {\Bbb Z}_{\ge -1}.
$$
We have $p_{-1}(1)=0, p_{0}(1)=1$ and $p_k(1) \geq 2, \forall k \geq 1$.
Therefore, if $p_{-1}$ can be realized on a $4$-fold, there would be
amples of coindex $4$ with $|H|=\emptyset$! It would be very interesting 
to find Fano $4$-folds like this, or to show that they cannot exist!
\endproclaim

\head 4. Effective Nonvanishing on Surfaces
\endhead
\proclaim{The Nonvanishing Problem} Let $(X,B_X)$ be a projective 
pair with Kawamata log terminal singularities, $D$ a Cartier nef 
divisor and $H$ a nef and big ${\Bbb R}$-divisor on $X$ such that
$$
D \equiv K_X+B_X+H.
$$
Is it true that $H^0(X,{\Cal O}_X(D)) \neq 0$? If the answer is 
affirmative, we say that nonvanishing holds.
\endproclaim
It follows from the Riemann-Roch formula that nonvanishing holds on 
curves. In the sequel we give 
sufficient conditions for nonvanishing to hold on surfaces. The main
ingredient is the Riemann-Roch formula again (i.e. explicit 
description of the coefficients), so this approach
does not generalize.   

\proclaim{Proposition 4.1} Let $(X,B_X)$ be a projective log-surface
with Kawamata log terminal singularities,
$D$ a Cartier nef divisor and $H$ a nef and big ${\Bbb R}$-divisor 
on $X$ such that
$
D \equiv K_X+B_X+H.
$

1) If $\chi({\Cal O}_X) \geq 0$, then $H^0(X,{\Cal O}_X(D)) \neq 0$.

2) If $\chi({\Cal O}_X) < 0$, then $X$ is birationally ruled and let $F$
be the general member of the ruling. Then $H^0(X,{\Cal O}_X(D)) \neq 0$
if one of the following holds:

2a) $H F >1$, or

2b) $H\equiv cD$ and $c>\frac{1}{2}$.
\endproclaim

\proclaim{Remark} 
The assumption on singularities is used for Kawamata-Viehweg vanishing only.
The case 2a) covers the absolute case, i.e. $B_X=0$, since
$H F=DF+2-B_X F>1$ if $B_X F < 1$. 
\endproclaim

\proclaim{Lemma 4.2} Let $\pi:S \to C$ be a fibration from a nonsingular surface
$S$ to the nonsingular curve $S$ such that the general fiber $F$ is a 
smooth rational curve. Let $H$ be a nef ${\Bbb R}$-divisor on $S$ such 
that $HF=1$. Then
$$
\chi({\Cal O}_S) \geq -\frac{1}{2}H(H+ K_S).
$$
and equality holds if $\pi$ is a ${\Bbb P}^1$-bundle.
\endproclaim

\demo{Proof of Proposition 4.1}

 By Kawamata-Viehweg vanishing, 
we have $H^j(X,{\Cal O}_X(D)) = 0 \ \forall j>0$,
so Riemann-Roch gives
$$
h^0(D)=dim H^0(X,{\Cal O}_X(D))=\frac{1}{2}D(H +B_X) + \chi({\Cal O}_X).
$$

1) If $D H=0$, then since $D$ is nef and $H$ is big, we have $D \equiv 0$, 
so $(X,B_X)$ is a log Fano and $h^0(D)=1$. 
Therefore we can assume $DH>0$, so 1) follows.

2) Assume now $\chi({\Cal O}_X) < 0$. Blowing-up $X$ we can assume from 
that $X$ admits a ${\Bbb P}^1$-fibration $\pi:X \to C$
over a curve.
\newline
2a) Let $F$ be the general fiber and $a=HF$, which is positive, 
since $H$ is big. We apply Lemma for $\frac{1}{a}H$ and get
$$
\align
&h^0(D) \geq \frac{1}{2}D(H +B_X)-\frac{1}{2a^2}(H^2+a H K_X) \\
&=\frac{1}{2a^2}(H+aD)((a-1)H+a B_X).
\endalign
$$
Therefore
$
h^0(D) \geq \frac{a-1}{2a^2}H^2 >0
$ if $HF>1$, so 2a) follows.
\newline
2b) Assume now that $DF=a>0$. Applying Lemma for $\frac{1}{a}D$ we obtain
$$
\align
&h^0(D) \ge \frac{1}{2}D(D-K_X)-\frac{1}{2a^2}(D^2+a D K_X) \\
&= \frac{a+1}{2a}D((1-\frac{1}{a})D-K_X) \\
&=  \frac{a+1}{2a}D(-\frac{1}{a}D+H+B_X). 
\endalign
$$
If $H \equiv cD, c>\frac{1}{2}$ then $-\frac{1}{a}D+H+B_X \equiv 
(c-\frac{1}{a})D+B_X$, so nonvanishing holds if $c-\frac{1}{a}>0$.
Otherwise, $a\le \frac{1}{c}<2$, so $a=DF=1$. 
\newline
We prove nonvanishing in 
this case by reducing the problem to a ${\Bbb P}^1$-bundle
over a curve. 
Let $\pi=\pi_0 \circ \mu$ be the factorization of $\pi$ as 
in Lemma 4.2. 

Step 1. There is a nef and big Cartier divisor $D_T$ on $T$ such that
$\mu^*(D_T)=D$ and $(1-c)D_T \equiv K_T+B_T$, where $(T,B_T)$ is
klt again (we just need the effectiveness of $B_T$). 

Indeed, after a finite number of Castelnouvo contractions,
we can assume that $D$ is positive on 
rational $(-1)$-curves of $X$. 
We claim that this implies $X=T$. If not, $\pi$ will have at least a singular
fiber, say $\pi^*(P)=\sum_{i=0}^s m_i E_i, P \in C$, where all $E_i$'s are
smooth rational negative curves, and one of them, $E_0$ say, is a $(-1)$-curve.
But $1=DF=\sum_{i=0}^s m_i D E_i$, and since $D E_0 >0$, we get
$m_0=1, D E_i=0, \forall i \ge 1$. Therefore $E_0$ is the only one 
$(-1)$-curve in $\pi^*(P)$, which implies $K_X E_i \ge 0 \forall i\ge 1$.
Hence $-2=K_X F=-1+\sum_{i=1}^s m_i K_X E_i \ge -1$ gives the 
contradiction. 

Step 2. Suffices to prove nonvanishing on $T$, so we assume 
$X=T$ and $D=D_T$. 
Let $D \equiv C_0+bF$, where we use the same notations as in
Lemma 4.2. We assume $h^0(D)=0$ in order to get a contradiction.
\newline
Since $DF=1$, Lemma 4.2 gives $h^0(D)=-DK_X$, i.e. $e-2b+2(g-1)=0$. 
Moreover, $\chi({\Cal O}_X(D))=\chi(\pi_*({\Cal O}_X(D)))=-e+b+2(1-g)$,
so $e-b+2(g-1)=0$. Therefore $b=0$. The inequalities $e \ge -g$ and 
$g \ge 2$ give us $g=2, e=-2$. 
\newline
Since $e=-2$, the cone of effective curves $NE(T)$ is included
in $\{aC_0+bF; a+b \geq 0\}$ \cite{Ha}. But $B_X \equiv (1-c)D-K_X \equiv
(3-c)C_0-4F$, which contradicts the effectiveness of $B_X$.

\hfill $\square$
\enddemo

\demo{Proof of Lemma 4.2}
From the classification of surfaces, it follows that 
$\pi=\pi_0 \circ \mu$,
where $\pi_0:T \to C$ is a ${\Bbb P}^1$-bundle and $\mu:S \to T$ is
a composition of blow-ups in nonsigular points. We prove the inequality
by induction on $\varrho(S/T)$, the number of blow-ups. 

If $\varrho(S/T)=0$, then equality holds.
Indeed, let $C_0$ be the minimal section, $C_0^2=-e, g=g(C)$ and let
$H \equiv C_0+bF$. The canonical divisor is 
 $K_S \equiv -2C_0+(2g-2-e)F$, hence
$H K_S=2(g-1)+e-2b$ and $H^2=-e+2b$. Therefore 
$H^2+H K_S=2(g-1)=-2 \chi({\Cal O}_S)$.

 Assume now $\varrho(S/T)>0$
and let $\nu:S \to S_1$ be the last blow-up of $\mu$. Let $E$ be the 
corresponding $(-1)$-curve on $S$ and $\alpha=H E \geq 0$. Then
there is a nef ${\Bbb R}$-divisor $H_1$ on $S_1$ such that 
$H+\alpha E=\nu^*H_1$ (note that $H_1 F=H F =1$). Therefore 
$H K_S+H^2=(\nu^*H_1-\alpha E)(\nu^* K_{S_1}+E)+
(\nu^*H_1-\alpha E)^2=H_1 K_{S_1}+H_1^2+
\alpha(1-\alpha)$. But $E$ appears in a singular $\pi$-fiber, 
which is equivalent to $F$.
Since $H$ is nef, we get $H F \geq H E$, that is $1\geq \alpha$.

Thus $H K_S+H^2 \geq H_1 K_{S_1}+H_1^2$,
and the inductive step finishes the proof.  
\hfill $\square$
\enddemo

\head 5. Proof of The main theorem
\endhead
\proclaim{Lemma 5.1} 
Let $|H|$ be a linear system on a klt log 
variety $(X,B_X)$ such that $dim\vert H \vert >0$
and let $S \in |H|$ a general member. 
Then the log canonical threshold 
$\gamma=lct(S;X,B_X)$ is constant (maximal) for general $S$, 
and $0< \gamma \leq 1$.
\newline
1) If $\gamma \neq 1$, then $LCS(X,B_X+\gamma S) \subseteq 
 Bsl |H|$.
\newline
2) If $\gamma = 1$, then one of the following holds:

(2.a) There is an lc center of 
   $(X,B_X+S)$ included in $Bsl |H|$, or

(2.b) $(X,B_X+S)$ is pure log terminal
      and $|H|$ has no fixed components.

Moreover, in case (2.b), $S$ is a (possibly disconnected) variety  
and $(S,B_S)$ is a klt pair, where $B_S=B_X\vert_S$.
\endproclaim
\demo{Proof}
 Resolving the singularities of $X$ and the base locus of $|H|$, there 
is a log resolution $\mu : Y \to X$ such that

(i) $\mu^* (K_X+B_X)=K_Y+ \mu^{-1} B_X + \sum a_j F_j$, with 
$a_j \in \Bbb Q, a_j < 1$,

(ii) $\mu^*|H|=|L|+\sum r_j F_j,$ with $|L|$ base point free, 
$r_j \in \Bbb N$ and 
       $r_j \neq 0$ iff $\mu(F_j) \subset Bsl|H|.$

Take $S \in |H|$ a general member. Then $\mu^{-1} S=
 T+ \underset{F_j fixed \ cpnt}\to {\sum} r_j F_j$, and $T$ is 
 smooth (possibly disconnected) from the first theorem of Bertini. Therefore
$$
\align
&\mu^*(K_X+ B_X+\gamma S)= K_Y +  \mu^{-1} B_X + 
  \gamma \mu^{-1} S \\
&+  \sum_{F_j exc, \mu F_j \subseteq Bsl}
 (a_j+ \gamma r_j) F_j
  + \sum_{F_j exc, \mu F_j \not \subseteq Bsl} a_j F_j.
\endalign  
$$
It is clear that places of log canonicity for $(X, B_X+\gamma S)$ 
can be just $T$, 
fixed components, or exceptional $F_j$'s over the base locus of $|H|$.

1) If $\gamma <1$, then $T$ cannot be an lc place, therefore
 all lc places lie over $Bsl|H|$.

2) If $\gamma=1$, $S$ is reduced and any fixed component would have 
 multiplicity $1$. Let us assume that there are no lc centers 
 included in $Bsl|H|$. Then lc places are exactly the components of $T$,
 hence there are no fixed components. Finally, the strict transform 
 $\mu^{-1}_{*} S=T$ is smooth, hence a disjoint union of components.
 Therefore (cf. \cite{Sh3}, pg. $99$) $(X,S)$ is plt.

Moreover, in case (2.b), every connected component of
$S$ is irreducible reduced,
otherwise we find $Z$ in $CLC(X,S)$ of codimension 
at least two at the intersection of two components, contradiction.
 Because $S$ is the unique center of log canonicity for $(X,B_X+S)$, it is 
minimal, hence normal \cite{Ka1}. Finally, $(X,B_X+S)$ implies
that $(S,B_S)$ is klt by the standard 
argument.
\hfill $\square$
\enddemo

\proclaim{Lemma 5.2} 
Let $(X,B_X)$ be a weak log Fano and let $H$ a Cartier nef and big
divisor with $-(K_X+B_X) \equiv ind(H) H$. Assume $(X,B_X+G)$ is 
log canonical with 
$LCS=LCS(X,B_X+G) \neq \emptyset$, and  $G \equiv \gamma H$. Then for 
$t>-ind(H) +\gamma$ we have:
\newline
i)$H^j(LCS,tH)=0, \ \forall j>0$,
\newline
ii) $H^0(X,tH) \longrightarrow H^0(LCS,tH|_{LCS}) \longrightarrow 0.$
\endproclaim

\demo{Proof}
 We can write $tH \equiv K_X+B_X+G + (t+ind(H) -\gamma)H$, 
so Shokurov's Theorem and its corollary (Appendix) applies for 
$t>- ind(H) +\gamma$. For the same values, 
 Kawamata-Viehweg vanishing gives $H^i(X,tH)=0 \ i>0$, hence the lemma 
follows.
\hfill $\square$
\enddemo
\proclaim{Lemma 5.3} In the same hypothesis as the previous lemma, assume that
$H$ is ample and there is a minimal log canonical center $Z$ of $(X,B_X+G)$ 
which is included in $Bsl|H|$. Then

1) Assume $H^0(Z,{\Cal O}_X(H))\ne 0$ for any ample Cartier divisor $H$ of 
coindex 
 $r-(codZ-\gamma)+\epsilon$($0< \epsilon \ll 1$) on  
 a log Fano of dimension $dimZ$. Then
 $$\gamma \ge ind(H)=n-r+1.$$

2) If $dimZ \leq 1$, then $\gamma \geq ind(H) +1$.

3)(cf. ~\cite{Me1} for $n=r=3$) We have $\gamma \geq codZ-r+2$.
\endproclaim
\demo{Proof}
 We show that $Z$ is not in the base locus if the opposite
(strict) inequalities hold. Perturbing 
$G$, we can assume that $CLC(X,B_X+G)=\{Z\}$ and that adjunction holds.
In particular, $LCS(X,B_X+G)=Z$, so Lemma 5.2 gives 
$$
H^0(X,{\Cal O}_X(H)) \longrightarrow H^0(Z,{\Cal O}_Z(H))\longrightarrow 0
\text{   if  } \gamma < ind(H)+1,
$$
therefore is enough to prove $H^0(Z,{\Cal O}_Z(H)) \neq 0$. By 
adjunction formula, there is a boundary $B_Z$ which makes $(Z,B_Z)$
a klt pair such that
$$
H\vert_Z \equiv (K_X+B_X+(ind(H)+1)H)\vert_Z \equiv   
  K_Z+B_Z+(ind(H)+1-\gamma)H\vert_Z.
$$
Note that $(ind(H)+1-\gamma)H\vert_Z$ is ample on $Z$
if $\gamma<ind(H)+1$.

1) If $\gamma < ind(H)$, then $(Z,B_Z)$ becomes a log Fano and $H$ has
coindex $co(H|_Z)=r-(codZ-\gamma)\le r+\epsilon$, so we have nonvanishing by
assumption.

2) Assume $\gamma<ind(H)+1$.
If $Z$ is a point, then $H^0(Z,{\Cal O}_Z(H))=k$, hence we're done. If
$dim(Z)=1$, then we know that nonvanishing holds for curves (Section 4), 
so we are done again.

3)(This argument works for $H$ nef and big)
Let $p(t)=\chi({\Cal O}_Z (tH))$, which is a polynomial of degree 
at most $dimZ$, and assume $dimZ> 0$. If $\gamma<codZ-r+2$, Lemma 5.2 
gives 
$p(t)=h^0(Z,tH)$ for $t\geq -(dimZ-1)$. In particular, $p(0)=1$.

If $deg p(t) <dimZ$, then, since it has at least $d-1$ zeros, $p(t) \equiv 0$.
Contradiction with $p(0)=1$.

If $deg p(t)=dimZ$, then $p(t)$ has $d-1$ negative zeros $-1,-2,\ldots,-(d-1)$, 
and $p(0)= 1$, hence $p(1)>0$.

\hfill $\square$
\enddemo

\demo{Proof of the Main Theorem}
 Since $r<4$, we have $dim|H|>0$ from Section 3 and let $S \in |H|$ 
be a general element, and $\gamma =lct(S;X,B_X), 0<\gamma \leq 1$.  
By Lemma 5.1, it is enough to prove that $(X,B_X+\gamma S)$
has no centers of log canonicity included in $Bsl|H|$.

Let us assume $Z$ is such a (minimal) lc center included in $Bsl|H|$
and derive a contradiction. We can also assume 
that $H$ ample by the Remark 2.6.

Since $\gamma<ind(H)+1$, the exact restriction sequence and the
adjunction formula from the proof of Lemma 5.3 hold true.
From Section 3 again, the assumptions of 5.3.1) are satisfied, so
$$
ind(H) \leq \gamma.
$$
This implies $n-r+1\le \gamma \le 1$, hence $n \le r< 4$.
Therefore $dim(X) \leq 3$, so $dim(Z)\leq 2$. By 5.3.2),
this implies $dim(X)=3, dim(Z)=2$.
But then $Z$ will be a fixed component of $|H|$ and
$\gamma=\frac{1}{m}$, where $m \in {\Bbb Z}_{\ge 1}$ is the 
multiplicity of $Z$ in $|H|$.
\newline
If $m\ge 2$, then $ind(H)+1-\gamma>\frac{1}{2}$, so Proposition 4.1.2b)
gives $H^0(X,{\Cal O}_X(H)) \neq 0$. Contradiction.
\newline
If $m=1$, then $(X,B_X+S)$ is exactly log canonical. Moreover,
$S \neq Z$ since $dim|H|>0$, and $S$ is connected since $H$ is ample.
Therefore $Z$ intersects another component of $S$, which
contradicts the minimality of the lc center $Z$.

 The existence of the ladder is clear. It is regular because
$-(K_{X_j}+B_{X_j})$ ($j>c$) is nef and big, so 
$H^1(X_j,{\Cal O}_{X_j})=0$ by Kawamata-Viehweg vanishing. 
Note also that $B_{X_{j-1}}=B_{X_j}|_{X_{j-1}}$ is well defined:
$B_{X_j}$ and $X_{j-1}$ have no common components, 
since $(X_j,B_{X_j}+X_{j-1})$ is plt.
\hfill $\square$
\enddemo

\demo{Proof of Corollary}
 The only nontrivial part is the preservation of smoothness.
Let $X > S > \cdots > X_2$ be the Fano ladder and assume $P$ is 
an isolated singularity of $X_2$. Let $E$ be the exceptional 
divisor on the blow-up of $X$ in $P$. Then
$$
a(P,X)=a(P,X_2)+ \sum_{i=2}^{n-1} m(P;X_i)> (n-2)m(P;S),
$$ 
where $a(P,X)$ denotes the log discrepancy of $X$ in $P$ and 
$m(P;S)$ is the multiplicity of $S$ in $P$.
But $a(P;X)=n$, hence $m(P;S) < 1+\frac{2}{n-2}$. If $n\geq 4$, we
get $m(P;S)=1$, so we are done. If $n=3$, the only nonsmooth case
is when $a(P;S)=1, m(P;S)=2$, i.e $P \in S$ is a canonical double point.
\hfill $\square$
\enddemo

\head 6. Appendix
\endhead

The following is an exposition of an unpublished result of V. Shokurov.
Let $(X,B_X)$ be a log pair and denote 
${\Cal I}={\Cal I}(X,B_X)$ and ${\Cal L}(X,B_X)$ 
the log canonical singularity subscheme with support 
$LCS(X,B_X)$. Note that ${\Cal L}(X,B_X)=LCS(X,B_X)$
if $(X,B_X)$ is log canonical. Let $D$ be a Cartier divisor on $X$.
 We have the following two exact sequences:
$$0 \longrightarrow {\Cal I} \longrightarrow {\Cal O}_X
            \longrightarrow {\Cal O}_{\Cal L(X,B_X)}  
      \longrightarrow 0,$$
$$0 \longrightarrow {\Cal I}(D) \longrightarrow {\Cal O}_X(D)
            \longrightarrow {\Cal O}_{\Cal L(X,B_X)}(D)  
 \longrightarrow 0.
$$
\proclaim{Theorem}(V. Shokurov)  
Let $(X,B_X)$ be a log pair, and let $D$ be a 
Cartier divisor and $H$ a nef and big $\Bbb R$-divisor such that
$$
D \equiv K_X+B_X+H.
$$
Then $H^j(X,{\Cal I}(D))=0, \forall j>0.$
\endproclaim
\demo{Proof}
 Let $f:(Y,B^Y) \to (X,B_X)$ be a log resolution. 
Then $f^*D \equiv K_Y+B^Y+f^*H$ and $f^*H$ is nef and big on $Y$. 
We use now the Grothendieck spectral sequence
associated to $Y \overset{f}\to{\to} X \overset{\pi}\to{\to} pt$:
$$
R^p \pi_*(R^q f_* {\Cal F}) \Longrightarrow R^{p+q}(\pi f)_*{\Cal F}, 
$$
for ${\Cal F}=f^*D+\lceil -B^Y \rceil$.  
By Kawamata-Viehweg vanishing, 
$R^i f_*({\Cal F})=0\forall i>0$, so the spectral sequence degenerates.
 Moreover,
 $R^0 f_*({\Cal F})=f_*(\lceil -B^Y \rceil) \otimes 
 {\Cal O}_X(D)={\Cal I}(D)$. Therefore
$$
H^i(X,{\Cal I}(D)) \simeq H^i(Y,{\Cal F}) \forall i,
$$
while $H^i(Y,{\Cal F})=0\forall i>0$ by Kawamata-Viehweg vanishing again.
\hfill $\square$
\enddemo
The short exact sequence gives the following extension of Kawamata-Viehweg
vanishing:
\proclaim{Corollary} In the same hypothesis as above, assume $(X,B_X)$
is a log canonical pair, and let $Y=LCS(X,B_X)$. Then the canonical map
$$
H^0(X,{\Cal O}_X(D))
            \longrightarrow H^0(Y,{\Cal O}_Y(D))
$$ is 
surjective and
$$
H^i(X,{\Cal O}_X(D))
            \simeq H^i(Y,{\Cal O}_Y(D))
$$ 
is an isomorphism for any $i>0$.
\endproclaim

\enddocument